\documentclass[10pt]{article}
\usepackage{amsmath,amsthm,amsfonts,amssymb}
\usepackage{epsfig}
\usepackage{color}
\usepackage{latexsym}
\usepackage{supertabular}
\usepackage{graphicx}
\usepackage{a4wide}
\numberwithin{equation}{section}

\def\whitebox{{\hbox{\hskip 1pt
 \vrule height 6pt depth 1.5pt
 \lower 1.5pt\vbox to 7.5pt{\hrule width
    3.2pt\vfill\hrule width 3.2pt}%
 \vrule height 6pt depth 1.5pt
 \hskip 1pt } }}
\def\qed{\ifhmode\allowbreak\else\nobreak\fi\hfill\quad\nobreak
     \whitebox\medbreak}

\newcommand{\ignore}[1]{}

\theoremstyle{plain}
\newtheorem{theorem}{Theorem}[section]
\newtheorem{corollary}[theorem]{Corollary}
\newtheorem{lemma}[theorem]{Lemma}

\newtheorem{remark}[theorem]{Remark}
\newtheorem{conjecture}[theorem]{Conjecture}
\def\qed{{\hfill$\square$}}
\def\proof{{\vspace{-0.3cm}\bf Proof: \,}}

\def\Z{{\mathbb Z}}
\def\Q{{\mathbb Q}}

\def\C{{\mathbb C}}
\def\F{{\mathbb F}}
\def\mod{{\mathrm{mod\,\,}}}

\def\Tr{{\mathrm{Tr}}}
\def\Norm{{\mathrm{Norm}}}
\def\Cay{{\mathrm{Cay}}}

\title{Lifting Constructions of Strongly Regular Cayley Graphs}
\author{Koji Momihara\footnotemark[1] \quad Qing Xiang\footnotemark[2]
}
\date{} 
\begin{document}
\maketitle
\renewcommand{\thefootnote}{\fnsymbol{footnote}}
\footnotetext[1]{
Department of Mathematics, Faculty of Education, Kumamoto University,  
2-40-1 Kurokami, Kumamoto 860-8555, Japan; Email: 
momihara@educ.kumamoto-u.ac.jp}
\footnotetext[2]{
Department of Mathematical Science, University of Delaware, Newark,DE 19716, USA; Email: xiang@math.udel.edu
}
\renewcommand{\thefootnote}{\arabic{footnote}}
\begin{abstract}  
We give two ``lifting'' constructions of strongly regular Cayley graphs. 
In the first construction we ``lift" a cyclotomic strongly regular 
graph by using a subdifference set of the Singer difference set. 
The second construction uses quadratic forms over finite fields and it 
is a common generalization of the construction of the affine polar graphs 
\cite{CK86} and a construction of strongly regular Cayley graphs given in 
\cite{FWXY}. The two constructions are related in the following way: 
The second construction can be viewed as a recursive construction, and the strongly regular Cayley graphs  obtained from the first construction can serve as starters for the second construction.  We also obtain association schemes from the second construction. 
\end{abstract}
\begin{center} 
\end{center}
\section{Introduction}
In this paper, we assume that the reader is familiar with the theory of strongly regular graphs and difference sets. For the theory of strongly regular graphs, our main references are \cite{bh} and \cite{gr}. For the theory of difference sets, we refer the reader to Chapter 6 of \cite{bjl}. Strongly regular graphs (srgs) are closely related to other combinatorial objects, such as two-weight codes, two-intersection sets in finite geometry, 
and partial difference sets. For these connections, we refer the reader to \cite[p.~132]{bh} and \cite{CK86, M94}.  

Let $\Gamma$ be a simple and undirected graph and $A$ be its adjacency matrix. A very useful way to check whether $\Gamma$ is strongly regular is by using the eigenvalues of $A$ (which are usually called eigenvalues of $\Gamma$). For convenience, we will call an eigenvalue of $\Gamma$ {\it restricted} if it has an eigenvector perpendicular to the all-ones vector ${\bf 1}$. Note that for a $k$-regular connected graph, the restricted eigenvalues are simply the eigenvalues different from $k$.

\begin{theorem}\label{char}
For a simple $k$-regular graph $\Gamma$ of order $v$, not complete or edgeless, with adjacency matrix $A$, the following are equivalent:
\begin{enumerate}
\item $\Gamma$ is strongly regular with parameters $(v, k, \lambda, \mu)$ for certain integers $\lambda, \mu$,
\item $A^2 =(\lambda-\mu)A+(k-\mu) I+\mu J$ for certain real numbers $\lambda, \mu$, where $I, J$ are the identity matrix and the all-ones matrix, respectively, 
\item $A$ has precisely two distinct restricted eigenvalues.
\end{enumerate}
\end{theorem}
One of the most effective methods for constructing srgs is by the Cayley graph construction. For example, the Paley graph ${\rm P}(q)$ is a class of well-known Cayley graphs on the finite field $\F_q$; that is, the vertices of ${\rm P}(q)$ are the elements of $\F_q$, and two vertices are adjacent if and only if their difference is a nonzero square. The parameters of ${\rm P}(q)$ are $(v,k,\lambda,\mu)=(4t+1,2t,t-1,t)$, where $q=4t+1$ is a prime power. More generally, let $G$ be an additively written group of order $v$, and let $D$ be a subset of $G$ such that $0\not\in D$ and $-D=D$, where $-D=\{-d\mid d\in D\}$. The {\it Cayley graph on $G$ with connection set $D$}, denoted by ${\rm Cay}(G,D)$, is the graph with the elements of $G$ as vertices; two vertices are adjacent if and only if their difference belongs to $D$. In the case where $\Cay(G,D)$ is strongly regular, the connection set $D$ is called a (regular) {\it partial difference set}.  The survey of Ma~\cite{M94} contains much of what is known about partial difference sets and about connections with strongly regular Cayley graphs. 

A classical method for constructing strongly regular Cayley graphs on the additive groups of finite fields is to use cyclotomic classes of finite fields. 
Let $p$ be a prime,  $f$ a positive integer, and let $q=p^f$. Let $e>1$ be an integer such that $e|(q-1)$, and $\gamma$ be a primitive element of $\F_q$. 
Then the cosets $C_i^{(e,q)}=\gamma^i \langle \gamma^e\rangle$, $0\leq i\leq e-1$, are called the {\it cyclotomic classes of order $e$} of $\F_q$. 
Many authors have studied the problem of determining when a union $D$ of cyclotomic classes forms a partial difference set. We call $\Cay(\F_q,D)$ a {\it cyclotomic strongly regular graph} if $D$ is a single cyclotomic class of $\F_q$ and $\Cay(\F_q,D)$ is strongly regular. Extensive work has been done on cyclotomic srgs, see \cite{BMW82,bwx,FMX11,FX111,GXY11,DL95, VLSch,M75,Momi,SW02,S67}. (Some of these authors used the language of cyclic codes in their investigations instead of strongly regular Cayley graphs or partial difference sets. We choose to use the language of srgs here.) The Paley graphs are primary examples of  cyclotomic srgs.  Also, if $D$ is the multiplicative group of a subfield of $\F_q$, then it is clear that $\Cay(\F_q , D)$ is strongly regular.  These cyclotomic srgs are usually called {\it subfield examples}. Next, if there exists a positive integer $t$ such that $p^t\equiv -1\,(\mod{e})$, then $\Cay(\F_q , D)$ is strongly regular. See \cite{BMW82}. These examples are usually called {\it semi-primitive}.  Schmidt and White made the following conjecture on cyclotomic srgs. 

\begin{conjecture}\label{con:SW}{\em (\cite{SW02})}
Let $\F_{p^f}$ be the finite field of order $p^f$, $e\,|\,\frac{p^f-1}{p-1}$ with $e>1$, and $C_0=C_{0}^{(e,p^f)}$ with $-C_0=C_0$. If $\Cay(\F_{p^f},C_0)$ is strongly regular, then one of the following holds: 
\begin{enumerate}
\item[(1)] (subfield case) $C_0=\F_{p^d}^\ast$ where $d\,|\,f$,  
\item[(2)] (semi-primitive case) $-1\in \langle p\rangle\le (\Z/e\Z)^\ast$,
\item[(3)] (exceptional case) $\Cay(\F_{p^f},C_0)$ has one of the eleven sets of parameters 
given in Table~\ref{Tab1}. 
\begin{table}[h]
\caption{Eleven sporadic examples}
\label{Tab1}
$$
\begin{array}{|c||c|c|c|c|}
\hline
\mbox{No.}&e&p&f&[(\Z/e\Z)^\ast:\langle p\rangle]\\
\hline
1&11&3&5&2\\
2&19&5&9&2\\
3&35&3&12&2\\
4&37&7&9&4\\
5&43&11&7&6\\
6&67&17&33&2\\
7&107&3&53&2\\
8&133&5&18&6\\
9&163&41&81&2\\
10&323&3&144&2\\
11&499&5&249&2\\
\hline
\end{array}
$$
\end{table}
\end{enumerate}
\end{conjecture}
A strongly regular graph is said to be of {\it Latin square type} (respectively, {\it negative Latin square type}) if $(v,k,\lambda,\mu)=
(n^2,r(n-\epsilon),\epsilon n+r^2-3\epsilon r,r^2-\epsilon r)$ and $\epsilon=1$ (respectively, $\epsilon=-1$). Typical examples of srgs of Latin square type or negative Latin square type come from nonsingular quadrics in the projective space ${\rm PG}(m-1, q)$, where $m$ is even. It seems that we know more examples of srgs of Latin square type than srgs of negative Latin square type, see \cite{dx}.  Our first main result in this paper is a construction of negative Latin square type strongly regular Cayley graphs $\Cay(\F_{q^2},D)$ by lifting a cyclotomic strongly regular graph on $\F_q$. The proof relies on the Davenport-Hasse lifting formula on Gauss sums. 

In our second main theorem, we will give a recursive construction of strongly regular Cayley graphs by using quadratic forms over finite fields under the assumption that certain strongly regular Cayley graphs exist on a ground field. This construction generalizes the following two constructions.
 
\begin{theorem}{\em (\cite{CK86})}\label{thm:affine}
Let $Q:V=\F_q^n\to \F_q$ be a nonsingular quadratic form, where $n$ is even and $q$ is an odd prime power,  and let $D=\{x\in \F_q^n\,|\,
Q(x)\mbox{ is a nonzero square in $\F_q$}\}$. Then, $\Cay(V,D)$ is a strongly regular graph (which is the so-called affine polar graph). 
\end{theorem}  
Feng et.al \cite{FWXY} gave the following construction using uniform cyclotomy. 

\begin{theorem}{\em (\cite{FWXY})}\label{thm:semipri}
Let $p$ be a prime, $e>2$, $q=p^{2jr}$, where $r\ge 1$, $e\,|\,(p^j+1)$, 
and $j$ is the smallest such positive integer. 
Let $Q:V=\F_q^n\to \F_q$ be a nonsingular quadratic form, where $n$ is even, and let $D_{C_i^{(e,q)}}=\{x\in \F_q^n\,|\,
Q(x)\in C_i^{(e,q)}\}$ for $0\le i\le e-1$. Then, $\Cay(V,D_{C_i^{(e,q)}})$ is strongly 
regular for all $0\le i\le e-1$.
\end{theorem}
  
The strongly regular Cayley graphs obtained in Section~\ref{fromSubDF} can be used as starters for the second construction. In this way, we obtain a few infinite families of strongly regualr Cayley graphs with Latin square type or negative Latin square type parameters. Furthermore, we discuss association schemes related to the second construction and obtain several new association schemes. 
\section{Background on Gauss sums and strongly regular Cayley graphs}

Let $p$ be a prime, $f$ a positive integer, and $q=p^f$. The canonical additive character $\psi$ of $\F_q$ is defined by 
$$\psi\colon\F_q\to \C^{\ast},\qquad\psi(x)=\zeta_p^{\Tr _{q/p}(x)},$$
where $\zeta_p={\rm exp}(\frac {2\pi i}{p})$ is a complex primitive $p$-th root of unity and $\Tr _{q/p}$ is the trace from $\F_q$ to $\F_p$. All complex characters  of $(\F_q,+)$ are given by $\psi_a$, where $a\in \F_q$. Here $\psi_a$ is defined by
\begin{equation}\label{additive}
\psi_a(x)=\psi(ax), \;\forall x\in \F_q. 
\end{equation}
For a multiplicative character 
$\chi_e$ of order $e$ of $\F_q$, we define the {\it Gauss sum} 
\[
G_f(\chi_e)=\sum_{x\in \F_q^\ast}\chi_e(x)\psi(x).
\] 
From the definition we see clearly that $G_f(\chi_e)\in\Z[\zeta_{ep}]$, the ring of algebraic integers in the cyclotomic field $\Q(\zeta_{ep})$. Let $\sigma_{a,b}$ be the automorphism of $\Q(\zeta_{ep})$ defined by 
\[
\sigma_{a,b}(\zeta_e)=\zeta_{e}^a, \qquad
\sigma_{a,b}(\zeta_p)=\zeta_{p}^b,
\]
where $\gcd{(a,e)}=\gcd{(b,p)}=1$. 
Below we list several basic properties of Gauss sums \cite{BEW97}: 
\begin{itemize}
\item[(i)] $G_f(\chi_e)\overline{G_f(\chi_e)}=q$ if $\chi_e$ is nontrivial;
\item[(ii)] $G_f(\chi_e^p)=G_f(\chi_e)$, where $p$ is the characteristic of $\F_q$; 
\item[(iii)] $G_f(\chi_e^{-1})=\chi_e(-1)\overline{G_f(\chi_e)}$;
\item[(iv)] $G_f(\chi_e)=-1$ if $\chi_e$ is trivial;
\item[(v)] $\sigma_{a,b}(G_f(\chi_e))=\chi_e^{-a}(b)G_f(\chi_e^a)$.
\end{itemize}
In general, explicit evaluations of Gauss sums are very difficult. There are only a few cases where the Gauss sums have been evaluated. 
The most well-known case is the {\it quadratic} case, i.e., the order of $\chi_e$ is two. 
The next simple case is the so-called {\it semi-primitive case} (also known as {\it uniform cyclotomy} or {\it pure Gauss sum}), where there exists an integer $j$ such that $p^j\equiv -1\,(\mod{e})$, where $e$ is the order of the multiplicative character involved. The explicit evaluations of Gauss sums in these cases are given in \cite{BEW97}. The next interesting case is the index $2$ case where the subgroup $\langle p\rangle$ generated by $p\in (\Z/{e}\Z)^\ast$ is of index $2$ in $(\Z/{e}\Z)^\ast$ and $-1\not\in \langle p\rangle $. In this case, 
it is known that $e$ can have at most two odd prime divisors. 
Many authors have investigated this case, see \cite{YX10} for a complete solution  to the problem of evaluating index $2$ Gauss sums. 
Recently, these index $2$ Gauss sums were used in the construction of new infinite families of strongly regular graphs. See \cite{FMX11,FX111}. 

Now we recall the following well-known lemma in algebraic graph theory (see e.g., \cite{bh}). 

\begin{lemma}\label{Sec3Le2}
Let $(G, +)$ be an abelian group and $D$ a subset of $G$ such that $0\not\in D$ and $D=-D$.  Then, the restricted eigenvalues of $\Cay(G,D)$ are given by $\psi(D)$, $\psi\in \widehat{G}\setminus \{\psi_0\}$,  
where $\widehat{G}$ is the character group of $G$ and $\psi_0$ is the trivial character. 
\end{lemma}

Let $q$ be a prime power and let $C_i^{(e,q)}=\gamma^i \langle \gamma^e\rangle$, $0\le i\le e-1$, be 
the cyclotomic classes of order $e$ of $\F_q$, where $\gamma$ is a fixed primitive root of $\F_q$.  In order to check whether a candidate subset $D=\bigcup_{i\in I}C_i^{(e,q)}$ is a connection set of a strongly regular Cayley graph, by Theorem~\ref{char} and Lemma~\ref{Sec3Le2}, it is enough to show that $\psi(aD)=\sum_{x\in D}\psi(ax)$, $a\in \F_q^\ast$, take exactly two values,  where $\psi$ is the canonical additive character of $\F_q$. Note that the sum $\psi(aD)$ can be expressed as a linear combination of 
Gauss sums  (cf. \cite{FX111}) by using the orthogonality of characters: 
\begin{eqnarray}\label{eigen}
\psi(aD)=\frac{1}{e}
\sum_{\chi\in C_0^{\perp}}G_f(\chi^{-1})
\sum_{i\in I}\chi(a\gamma^i ), 
\end{eqnarray}
where 
$C_0^{\perp}$ is the subgroup of $\widehat{\F_q^\ast}$
consisting of all $\chi$ which are trivial on $C_0^{(e,q)}$. 
Thus, the computations needed to show whether a candidate subset $D=\bigcup_{i\in I}C_i^{(e,q)}$ is a connection set of a strongly regular Cayley graph are essentially reduced to evaluating Gauss sums. 
However, as previously said, evaluating Gauss sums explicitly is very difficult. In Section 3 of this paper, we will give a contruction of strongly regular graphs by ``lifting" a cyclotomic srg. To prove that our construction indeed gives rise to srgs, we do not evaluate the Gauss sums involved; instead, we use the Davenport-Hasse lifting formula stated below.

\begin{theorem}\label{thm:lift}{\em (\cite{BEW97})}
Let $\chi$ be a nontrivial multiplicative character of $\F_q=\F_{p^f}$ and 
let $\chi'$ be the lift of $\chi$ to $\F_{q'}=\F_{p^{fs}}$, i.e., $\chi'(\alpha)=\chi(\Norm_{q'/q}(\alpha))$ for $\alpha\in \F_{q'}$, where $s\geq 2$ is an integer. Then 
\[
G_{fs}(\chi')=(-1)^{s-1}(G_f(\chi))^s. 
\]
\end{theorem}
\section{Lifting cyclotomic strongly regular graphs via subdifference sets}\label{fromSubDF}

In this section, we give a construction of strongly regular Cayley graphs by ``lifting" a cyclotomic strongly regular graph via a subdifference set of the Singer difference sets. We start by reviewing a construction of the Singer difference sets. 

Let $p$ be a prime, $f\geq 1$, $m\geq 3$ be integers and $q=p^f$. Let $L$ be a complete system of coset representatives of $\F_q^\ast$ in $\F_{q^m}^\ast$. We may assume that $L$ is chosen in such a way that $\Tr_{q^m/q}(x)=0$ or $1$ for any $x\in L$. 
Let 
\[
L_0=\{x\in L\,|\,\Tr_{q^m/q}(x)=0\}\mbox{ and } L_1=\{x\in L\,|\,\Tr_{q^m/q}(x)=1\}. 
\]
Then, 
\begin {equation}\label{Singer}
H_0=\{\overline{x}\in \F_{q^m}^\ast/\F_q^\ast\,|\,x\in L_0\}
\end{equation}
is a Singer difference set. 

Note that any nontrivial multiplicative character $\chi$ of exponent $(q^m-1)/(q-1)$ of $\F_{q^m}^\ast$ induces a character of the quotient group $\F_{q^m}^\ast/\F_q^\ast$, which will be denoted by $\chi$ also. Moreover, every character of $\F_{q^m}^\ast/\F_q^\ast$ arises in this way. By a result of Yamamoto \cite{Yamamoto}, for any nontrivial multiplicative character $\chi$ of exponent $(q^m-1)/(q-1)$ of $\F_{q^m}^\ast$, we have
$$\chi(H_0)=G_{fm}(\chi)/q.$$  

Now, let $\F_q^\ast\le C_0(:=C_0^{(e,q^m)})\le \F_{q^m}^\ast$ be a subgroup such that $[\F_{q^m}^\ast : C_0]=e$. Then $$\overline{C_0}=C_0/\F_q^\ast\le \F_{q^m}^\ast/\F_q^\ast.$$ 
Let $S$ be a complete system of coset representatives of $\overline{C_0}$ in $\F_{q^m}^\ast/\F_q^\ast$,  and $G=\{\overline{s}\,|\,s\in S\}\simeq \F_{q^m}^\ast/C_0$. 
Then, by assumption, $[\F_{q^m}^\ast : C_0]$ divides $(q^m-1)/(q-1)$, i.e., $e=|G|\mid (q^m-1)/(q-1)$. 

From now on, we assume that $\Cay(\F_{q^m},C_0) $ is strongly regular. Then $|H_0\cap s\overline{C_0}|$, $s\in S$, take exactly two values. (See  \cite{CK86} or \cite{SW02}.) It follows that $|H_0\cap s\overline{C_0}|-|H_0\cap \overline{C_0}|=0$ or $\delta$, where $\delta$ is a nonzero integer. For any nontrivial multiplicative character $\chi$ of $\F_{q^m}$ of exponent $e$, we have 
\begin{eqnarray*}
\chi(H_0)&=&\sum_{s\in S}|H_0\cap s\overline{C_0}|\chi(\overline{s})\\
&=&\sum_{s\in S}(|H_0\cap s\overline{C_0}|-|H_0\cap \overline{C_0}|)\chi(\overline{s})\\
&=&\delta\sum_{s\in S'}\chi(\overline{s}), 
\end{eqnarray*}
where 
\begin{equation}\label{def:S'}
S'=\{s\in S : |H_0\cap s\overline{C_0}|-|H_0\cap \overline{C_0}|=\delta\}.
\end{equation}
 Thus
\begin{equation}\label{eq:sum}
\sum_{s\in S'}\chi(\overline{s})=\frac{\chi(H_0)}{\delta}=\frac{G_{fm}(\chi)}{\delta q}. 
\end{equation}
It follows that $\delta$ is a power of $p$, and $\overline{S'}:=\{\overline{s}\mid s\in S'\}\subset G$ is a $(e, |S'|, \lambda')$-difference set, 
which is usually called a {\it subdifference set} of $H_0$. See Section 6 of \cite{SW02}. The term ``subdifference set'' was first introduced by McFarland \cite{Mc}. 

Let $\gamma$ be a primitive element of $\F_{q^{2m}}$ and let $\omega=\Norm_{q^{2m}/q^{m}}(\gamma)=\gamma^{q^m+1}$, which is a
primitive element of the subfield $\F_{q^m}$ of $\F_{q^{2m}}$. Let 
$C_j^{(e,q^{2m})}=\gamma^j\langle \gamma^e\rangle$ and  
$C_j^{(e,q^m)}=\omega^j\langle \omega^e\rangle=\omega^j C_0$. 

\begin{theorem} \label{theorem:main1}
Assume that $\F_q^\ast\le C_0\le \F_{q^m}^\ast$ be a subgroup such that $[\F_{q^m}^\ast : C_0]=e$, $-C_0=C_0$, and $\Cay(\F_{q^m},C_0)$ is strongly regular. Let 
$I=\{0\le i\le e-1\mid \overline{\omega}^i\in S'\}$, where $S'$ is defined in (\ref{def:S'}) and $\overline{\omega}$ stands for $\omega\F_q^{\ast}$. 
Let 
$$D=\bigcup_{i\in I}C_i^{(e,q^{2m})}.$$ 
Then $\Cay(\F_{q^{2m}},D)$ is also strongly regular.    
\end{theorem}

\proof
Let $\psi_1$ be the canonical additive character of $\F_{q^{2m}}$ and let  $\chi_e'$ be a multiplicative character of order $e$ of $\F_{q^{2m}}$. The restricted eigenvalues of $\Cay(\F_{q^{2m}},D)$ are $\psi_1(\gamma^aD)$, $0\leq a\leq e-1$. By (\ref{eigen}), in order to show that $\Cay(\F_{q^{2m}},D)$ is strongly regular, we compute the sums 
\begin{eqnarray*}
T_a=e\cdot \psi_1(\gamma^a D)+|I|=\sum_{x=1}^{e-1}G_{2fm}(\chi_e'^{-x})\sum_{i\in I}\chi_e'^{x}(\gamma^{a+i}),
\end{eqnarray*}
where $a=0,1,\ldots,e-1$. 
Since $e\,|\,(q^m-1)$, $\chi_e'$ must be the lift of a character, say $\chi_e$, of $\F_{q^m}$. 
By the Davenport-Hasse lifting formula, we have 
\begin{eqnarray*}
T_a
=-\sum_{x=1}^{e-1}{\chi_{e}}^{x}(\omega^{a})G_{fm}(\chi_{e}^{-x})G_{fm}(\chi_{e}^{-x})
\sum_{i\in I}\chi_{e}^{x}(\omega^{i})
\end{eqnarray*}
By the definition of $I$, we have 
\[
\sum_{i\in I}\chi_e^x(\omega^i)=\sum_{s\in S'}\chi_e^x(s)=\frac{G_{fm}(\chi_e^x)}{\delta q}. 
\]
Hence
\begin{eqnarray}\label{exp}
T_a
&=&-\frac{1}{\delta q} \sum_{x=1}^{e-1}\chi_{e}^{x}(\omega^{a})G_{fm}(\chi_{e}^{-x})G_{fm}(\chi_{e}^{-x})G_{fm}(\chi_{e}^{x})\nonumber\\
&=&- \frac{q^{m-1}}{\delta} \sum_{x=1}^{e-1}\chi_{e}^{x}(\omega^{a})G_{fm}(\chi_{e}^{-x}), 
\end{eqnarray}
where in the last step we used  the fact that $G_{fm}(\chi_{e}^{x})G_{fm}(\chi_{e}^{-x})
=\chi_{e}^{x}(-1)q^m$.  
By the assumption that $\Cay(\F_{q^m}, C_0^{(e,q^m)})$ is strongly regular, we have $\sum_{x=1}^{e-1}\chi_{e}^{x}(\omega^{a})G_{fm}(\chi_{e}^{-x})$, $a=0,1,\ldots,e-1$, take exactly two values. We conclude that $T_a$, $0\leq a\leq e-1$, take exactly two values too. Therefore $\Cay(\F_{q^{2m}}, D)$ is also strongly regular. \qed

Note that the set $D$ has size $|D|=\frac{(q^m-1)}{e}|I|(q^m+1)$. By applying Theorem~\ref{theorem:main1} to the known cyclotomic srgs in the statement of Conjecture~\ref{con:SW}, 
we obtain a lot of strongly regular Cayley graphs. We first apply  Theorem~\ref{theorem:main1} to the semi-primitive examples. In this case, 
we have $|I|=|S'|=1$ by \cite[p.~23]{Del}.     

\begin{corollary}\label{cor:semi}
Let $p$ be a prime, $e\ge 2$, $q^m=p^{2jr}$, where $m=2jr$, $r\ge 2$, $e\,|\,(p^j+1)$, 
and $j$ is the smallest such positive integer. Then there exists an $(n^2,r(n+1),-n+r^2+3r,r^2+r)$ strongly regular Cayley graph with $n=q^m$ and $r=(q^m-1)/e$.
\end{corollary}

The proof is straightforward. We omit it. Next we apply Theorem~\ref{theorem:main1} to the subfield examples. 

\begin{corollary}\label{cor:singer}
Let $q$ be a prime power, $m\ge 3$  a positive integer and $a$ any positive divisor of $m$. 
Then there exists an $(n^2,r(n+1),-n+r^2+3r,r^2+r)$ strongly regular Cayley graph with 
$n=q^m$ and $r=q^{m-a}-1$.
\end{corollary}

\proof 
We apply Theorem~\ref{theorem:main1} to the subfield examples of cyclotomic srgs. We use the same notation as in the statement and proof of Theorem~\ref{theorem:main1}. Then, by \cite[p.~23]{Del}, we have $C_0=\F_{q^a}^\ast$, $e=\frac{q^m-1}{q^a-1}$, $|I|=|S'|=\frac{q^{m-a}-1}{q^a-1}$, $\delta=q^{a-1}$ and the restricted eigenvalues of $\Cay(\F_{q^m},C_0)$ are $-1$ and $q^a-1$. The corollary now follows by straightforward computations using (\ref{exp}).
\qed

\begin{remark}
In the case where $q=2$ and $a=1$, the parameters of the strongly regular Cayley graphs obtained in Corollary~\ref{cor:singer} are 
\[
(2^{2m},(2^m+1)(2^{m-1}-1),2^{m-1}(2^{m-1}-1)-2,2^{m-1}(2^{m-1}-1)).
\] 
Then, the set $D\cup \{0\}$ clearly forms a difference set  with parameters $(2^{2m},2^{m-1}(2^m-1),2^{m-1}(2^{m-1}-1))$, which is a Hadamard difference set in the elementarty abelian 2-group of order $2^{2m}$. This difference set was first discovered in \cite[p.~105]{dillon}. The corresponding bent function is a monomial quadratic bent function.
\end{remark}

Finally, we apply Theorem~\ref{theorem:main1} to the eleven sporadic  examples of cyclotomic srgs. In this case,  the values of $k:=|S'|$ are 
given in \cite[Table II]{SW02}.  
\begin{corollary}\label{cor-spor}
There exists a $(q^{2}, r(q+1),\lambda,\mu)$ negative Latin square type strongly regular Cayley graph, where $r=k(q-1)/e$, in each of the following cases:  
\begin{eqnarray*}
(q,e,k)&=&(3^5,11,5),(5^9,19,9),(3^{12},35,17),(7^9,37,9),(11^7,43,21),(17^{33},67,33)\\
& &(3^{53},107,53),(5^{18},133,33),(41^{81},163,81),(3^{144},323,161),(5^{249},499,249). 
\end{eqnarray*}
\end{corollary}

\section{Strongly regular Cayley graphs from quadratic forms}
\label{sec4}
Let $V$ be an $n$-dimensional vector space over $\F_q$. A function $Q:V\to \F_q$ is called a {\it quadratic form} if\\
(i) $Q(\alpha v)=\alpha^2 Q(v) $ for all $\alpha\in \F_q$ and $v\in V$,\\ 
(ii) the function $B:V\times V\to \F_q$ defined by 
$B(u,v)=Q(u+v)-Q(u)-Q(v)$ is bilinear. 

We say that $Q$ is {\it nonsingular} if the subspace $W$ of $V$
with the property that $Q$ vanishes on $W$ and $B(v,w)=0$ for all $v\in V$ and $w\in W$ is the zero subspace (equivalently, 
we say that $Q$ is nonsingular if it can not be written as a form 
in fewer than $n$ variables by any invertible linear change of variables).  
If $\F_q$ has odd characteristic or $V$ is even-dimensional over an even-characteristic field $\F_q$, then $Q$ is nonsingular 
if and only if $B$ is nondegenerate \cite[p.~14]{C99}. But this is not necessarily true in general.  
Now assume that $n$ is even if $q$ is even and  $n$ is arbitrary otherwise.  
Then, $Q: V=\F_q^n\to \F_q$ is a nonsingular quadratic form if and only if the associated polar form $B(x,y)=Q(x+y)-Q(x)-Q(y)$ is nondegenerate; 
the characters $\phi_b$, $b\in V$,  of $(V, +)$ defined by 
\begin{equation}\label{defchar}
\phi_b(x)=\psi_1(B(b,x)),\;  \forall x\in V,
\end{equation}
where $\psi_1$ is the canonical additive character of $\F_q$, are all the complex characters of $(V,+)$. We can linearly extend the characters $\phi_b$ to the whole group ring $\C[V]$:  for $A=\sum_{g\in V}a_g g\in \C[V]$, we define $\phi_b(A)=\sum_{g\in V}a_g \phi_b(g)$. 

It is well known that a nonsingular quadratic form on $V=\F_q^n$, where $n$ is even, is equivalent to either 
\begin{equation}\label{hyper}
x_1x_2+x_3x_4+\cdots+x_{n-1}x_n,
\end{equation}
or
\begin{equation}\label{ellip}
x_1x_2+x_3x_4+\cdots+x_{n-3}x_{n-2}+(a x_{n-1}^2+bx_{n-1}x_n+cx_n^2),
\end{equation}
where $a x_{n-1}^2+bx_{n-1}x_n+cx_n^2$ is irreducible over $\F_q$. 

A nonsingular quadratic form equivalent to (\ref{hyper})  (resp. (\ref{ellip})) is called {\it  hyperbolic} (resp. {\it elliptic}).
 
\begin{lemma}{\em (\cite[Theorem 3.2]{LS99})\label{lemma:Leep}}
Let $q=p^f$, where $p$ a prime and $f\geq 1$ is an integer, and let $Q$ be a nonsingular quadratic form on $V=\F_q^n$ with $n=2m$ even. Then 
\[
\sum_{x\in V}\psi_1(Q(x))=\epsilon q^{m}, 
\]
where $\epsilon=1$ or $-1$ according as  $Q$ is hyperbolic or  elliptic. 
\end{lemma}

For each $u\in \F_q$, define $D_u=\{x\in V\,|\,Q(x)=u\}$, and we use the same $D_u$ to denote the corresponding group ring element $\sum_{z\in D_u}z\in \C[V]$.  For a subset $X$ of $\F_q$, we write 
$D_X=\sum_{x\in X}D_x$, which is viewed as an element of $\C[V]$. Now, we give the following key lemma. 
\begin{lemma}\label{essen}
Let $q=p^f$ be a prime power and $n=2m$ be an even positive integer. Let $Q: V=\F_q^n\rightarrow \F_q$ be a nonsingular quadratic form. For any $e\,|\,(q-1)$, let $C_i^{(e,q)}=\omega^i\langle\omega^e\rangle$ and $C_{i}^{(e,q^2)}=\gamma^i\langle \gamma^e\rangle$, $0\leq i\leq e-1$, denote the cyclotomic classes of order $e$ of $\F_q$ and $\F_{q^2}$, respectively, where $\gamma$ is a fixed primitive element of $\F_{q^2}$ and $\omega=\Norm_{q^2/q}(\gamma)$. Then, for any $0\neq b\in V$,
\[
\phi_b(D_{C_i^{(e,q)}})=\left\{
   \begin{array}{ll}
  -\epsilon q^{m-1} \frac{q-1}{e}, & \mbox{if $Q(b)=0$,} \\
  -\epsilon q^{m-1} \psi_1'(\gamma^{i+s}C_0^{(e,q^2)}), & \mbox{if $Q(b)\in C_s^{(e,q)}$ for $0\leq s\leq e-1$,}
 \end{array}
   \right.
\]
and 
\[
\phi_b(D_{0})=\left\{
   \begin{array}{ll}
  \epsilon q^{m-1} (q-1), & \mbox{if $Q(b)=0$,} \\
  -\epsilon q^{m-1}, & \mbox{if $Q(b)\not=0$,}
 \end{array}
   \right.
\]
where $\epsilon=1$ or $-1$ according as  $Q$ is hyperbolic or  elliptic, $\phi_b$ is defined in (\ref{defchar}), and $\psi_1'$ is the canonical additive character of $\F_{q^2}$. 
\end{lemma}
\proof 
We compute  
the values of $\phi_b(D_{C_i^{(e,q)}})$. For $b\in V\setminus\{0\}$, we have
\begin{eqnarray}\label{eqn:1}
q\cdot \phi_b(D_{C_i^{(e,q)}})&=&\sum_{y\in C_i^{(e,q)}}\sum_{x\in V}\chi_b(x)\sum_{u\in \F_q}\psi_1(u(Q(x)-y))\nonumber\\
&=&\sum_{x\in V}\sum_{u\in \F_q}\psi_1(B(b,x)+uQ(x))\psi_u(-C_i^{(e,q)})\nonumber\\
&=&\sum_{x\in V}\sum_{u\in \F_q^\ast}\psi_1(B(b,x)+uQ(x))\psi_u(-C_i^{(e,q)})+\frac{q-1}{e}\sum_{x\in V}\psi_1(B(b,x)).\nonumber
\end{eqnarray}
Since $\sum_{x\in V}\psi_1(B(b,x))=0$ and 
$B(b,x)+uQ(x)=-u^{-1}Q(b)+uQ(x+u^{-1}b)$ for $u\in \F_q^{\ast}$, we have 
\begin{equation}\label{eq:qqq}
q\cdot \phi_b(D_{C_i^{(e,q)}})=\sum_{x\in V}\sum_{u\in \F_q^\ast}\psi_1(-u^{-1}Q(b)+uQ(x+u^{-1}b))\psi_u(-C_i^{(e,q)}). 
\end{equation}
By Lemma~\ref{lemma:Leep}, we obtain 
\begin{eqnarray*} 
q\cdot \phi_b(D_{C_i^{(e,q)}})&=&\epsilon q^{m} \sum_{u\in \F_q^\ast}\psi_1(-u^{-1}Q(b))\psi_u(-C_i^{(e,q)})\\
&=&\epsilon q^{m} \sum_{a=0}^{e-1}\sum_{c=0}^{(q-1)/e-1}\psi_1(\omega^{-a-ce}Q(b))\psi_1(\omega^{a+ce} C_i^{(e,q)})\\
&=&\epsilon q^{m} \sum_{a=0}^{e-1}\psi_1(\omega^{a+i}C_0^{(e,q)})\sum_{c=0}^{(q-1)/e-1}\psi_1(\omega^{-a-ce}Q(b))\\
&=&\left\{
   \begin{array}{ll}
  -\epsilon q^{m} \frac{q-1}{e}, & \mbox{if $Q(b)=0$,} \\
  \epsilon q^{m} \sum_{a=0}^{e-1}\psi_1(\omega^{a+i}C_0^{(e,q)})\psi_1(\omega^{-a+s}C_0^{(e,q)}), & \mbox{if $Q(b)\in C_s^{(e,q)}$ for  $0\leq s\leq e-1$.}
 \end{array}
   \right.
\end{eqnarray*}
Below we further prove that
\[
\sum_{a=0}^{e-1}\psi_1(\omega^{a+i}C_0^{(e,q)})\psi_1(\omega^{-a+s}C_0^{(e,q)})=-\psi_1'(\gamma^{i+s}C_0^{(e,q^2)}). 
\]
Let $\chi_e'$ be a multiplicative character of order $e$ of $\F_{q^2}$. 
Since $e\,|\,(q-1)$, $\chi_e'$ must be the lift of a character, say $\chi_e$, of $\F_{q}$. 
Then, by the orthogonality of characters and the Davenport-Hasse lifting formula on Gauss sums, 
we have 
\begin{eqnarray*}
& &\sum_{a=0}^{e-1}\psi_1(\omega^{a+i}C_0^{(e,q)})\psi_1(\omega^{-a+s}C_0^{(e,q)})\\
&=&\frac{1}{e^2}\sum_{x=0}^{e-1}\sum_{y=0}^{e-1}
G_f(\chi_e^x)G_f(\chi_e^y)\chi_e^{-x}(\omega^i)\chi_e^{-y}(\omega^{s})\left(\sum_{a=0}^{e-1}\chi_e^{-x+y}(\omega^{a})\right)\\
&=&\frac{1}{e}\sum_{x=0}^{e-1}
G_f(\chi_e^x)^2\chi_e^{-x}(\omega^{i+s})\\
&=&-\frac{1}{e}\sum_{x=0}^{e-1}
G_{2f}({\chi_e'}^x){\chi_e'}^{-x}(\gamma^{i+s})\\
&=&-{\psi'}_1(\gamma^{i+s}C_0^{(e,q^2)}).
\end{eqnarray*}
Similarly, for $b\in V\setminus\{0\}$, we have 
\[
q\cdot \phi_b(D_{0})=\left\{
   \begin{array}{ll}
  \epsilon q^{m} (q-1), & \mbox{if $Q(b)=0$,} \\
  -\epsilon q^{m}, & \mbox{if $Q(b)\not=0$.}
 \end{array}
   \right.
\]
The proof is now complete. \qed

Now we give the main theorem of this section. 
\begin{theorem}\label{theorem:main2}
Let $q$ be a prime power, $e>1$ be an integer dividing $q-1$, and  $I$ be a subset of $\{0,1,\ldots,e-1\}$.
Let $\gamma$ be a fixed primitive element of $\F_{q^2}$ and $\omega=\Norm_{q^2/q}(\gamma)$. Let  $C_{i}^{(e,q^2)}=\gamma^i\langle \gamma^{e}\rangle$  and  $E=\bigcup_{i\in I}C_i^{(e,q)}$. 
Assume that $\Cay(\F_{q^2},\bigcup_{i\in I}C_i^{(e,q^2)})$ is a negative Latin square type srg. 
Then, for any nonsingular quadratic form $Q:V=\F_q^n \to \F_q$, where $n=2m$, the Cayley graph $\Cay(V,D_{E})$  
is strongly regular. 
\end{theorem}
\proof 
By assumption, $\Cay(\F_{q^2},\bigcup_{i\in I}C_i^{(e,q^2)})$ is an srg of negative Latin square type. Its parameters are $
(q^2,r(q+1),-q+r^2+3r,r^2+r)$, where $r=|I|(q-1)/e$. The restricted eigenvalues of this srg are $r$ and $r-q$. 
It follows from Lemma~\ref{essen} that $\phi_b(D_E)=\sum_{i\in I}\phi_b(D_{C_i^{(e,q)}})$, $b\in V\setminus \{0\}$, take exactly two values, 
namely $-\epsilon q^{m-1}|I|(q-1)/e$ and $-\epsilon q^{m-1}(|I|(q-1)/e-q)$. Thus $\Cay(V,D_{E})$ is strongly regular. \qed 
\begin{remark}\label{re:sec4}
\begin{itemize}
\item[(i)]
Note that the srg $\Cay(V,D_{E})$ obtained in the above theorem is of Latin square type or negative Latin square type according as $\epsilon=1$ or $-1$. One can see that under the same assumptions as in Theorem~\ref{theorem:main2}, $\Cay(V,D_{E}\cup (D_{0}\setminus \{0\}))$ is also strongly regular 
since 
$\Cay(\F_{q^2},\bigcup_{\overline{I}}C_i^{(e,q^2)})$ is a negative Latin square type srg, where $\overline{I}=\{0,1,\ldots,e-1\}\setminus I$. Furthermore, it is well known that $\Cay(V,D_{0})$ is also strongly 
regular with Latin square type or negative Latin square type parameters according to $\epsilon=1$ or $-1$\cite{CK86}. 

\item[(ii)] The most important condition in Theorem~\ref{theorem:main2} is that $\Cay(\F_{q^2},\bigcup_{i\in I}C_i^{(e,q^2)})$ is strongly regular. This condition is trivially satisfied in the following case. Let $e=2$ and $I=\{0\}$. Then, 
$\Cay(\F_{q^2},C_0^{(2,q^2)})$ is obviously strongly regular with negative Latin square parameters. Thus, the aforementioned condition is trivially satisfied. In this case, the srg $\Cay(V, D_E)$ obtained from Theorem~\ref{theorem:main2} is exactly the affine polar graph. 
We thus have recovered Theorem~\ref{thm:affine}. 
\end{itemize}
\end{remark}

The strongly regular Cayley graphs obtained in Section~\ref{fromSubDF} satisfy the assumptions of Theorem~\ref{theorem:main2}, namely, $e$ divides $q-1$ and $\Cay(\F_{q^2},\bigcup_{i\in I}C_i^{(e,q^2)})$ is a  negative Latin square type srg. Thus, we can use the srgs obtained in Section~\ref{fromSubDF} as starters to obtain new ones. We first 
consider the semi-primitive case. 
Let $p$ be a prime, $e>2$, $q=p^{2jr}$, where $r\ge 1$, $e\,|\,(p^j+1)$, and $j$ is the smallest such positive integer. In this case, $E$ is chosen as the dual of the semi-primitive cyclotomic strongly regular Cayley graph $\Cay(\F_q,C_0^{(e,q)})$. Then, $|E|=(q-1)/e$ by \cite[p.~23]{Del}. (Here, replace $q^m$ of Corollary~\ref{cor:semi} 
with $q$.) 
By applying Theorem~\ref{theorem:main2} to srgs of Corollary~\ref{cor:semi}, we have the 
following corollary. 
\begin{corollary}
Let $p$ be a prime, $e>2$, $q=p^{2jr}$, $r\ge 2$, $e\,|\,(p^j+1)$, 
and $j$ is the smallest such positive integer. Then there exists a $(q^{2m},r(q^m-\epsilon),\epsilon q+r^2-3\epsilon r,r^2-\epsilon r)$ strongly regular Cayley graph with $r=q^{m-1}(q-1)/e$.
\end{corollary}
We have thus recovered Theorem~\ref{thm:semipri}.
Next we consider the subfield case. Let $q=p^{st}$ and $e=\frac{p^{st}-1}{p^s-1}$. Here, $E$ is chosen as the dual of the subfield cyclotomic strongly regular Cayley graph $\Cay(\F_q,C_0^{(e,q)})$. Then, $|E|=p^{s(t-1)}-1$ by \cite[p.~23]{Del}. (We replaced $n=q^m$ and $r=q^{m-a}-1$ in Corollary~\ref{cor:singer} with
$n=p^{st}$ and $r=p^{s(t-1)}-1$ respectively.) 
Then, we have the 
following corollary.  
\begin{corollary}
There exists a $(p^{2stm},(p^{stm}-\epsilon)r,\epsilon p^{stm}+r^2-3\epsilon r,r^2-\epsilon r)$ strongly regular Cayley graphs for any prime $p$ and positive integers $s,t,$ and $m$, where $\epsilon=\pm 1$ and $r=p^{st(m-1)}(p^{s(t-1)}-1)$. 
\end{corollary} 
Finally, we consider the sporadic cases. Let $q,e$ and $k$ be as those in Corollary~\ref{cor-spor}. In this case,  $E$ is chosen as the dual of sporadic cyclotomic strongly regular Cayley graphs $\Cay(\F_q,C_0^{(e,q)})$. Then, $|E|=\frac{(q-1)}{e}\cdot k$. 
(Note that the number $k$ is the size of the subdifference set corresponding to the cyclotomic srg $\Cay(\F_{q},C_0^{(e,q)})$, see \cite[Table II]{SW02}.) 
Hence, we obtain the following corollary. 
\begin{corollary}
There exists a $(q^{2m}, r(q^m-\epsilon),\epsilon q^{m}+r^2-3\epsilon r,r^2-\epsilon r)$ strongly regular Cayley graph for any $m\ge 1$, 
where $\epsilon=\pm 1$ and  $r=q^{m-1}\cdot \frac{(q-1)}{e}\cdot k$, 
in each of the following cases:  
\begin{eqnarray*}
(q,e,k)&=&(3^5,11,5),(5^9,19,9),(3^{12},35,17),(7^9,37,9),(11^7,43,21),(17^{33},67,33)\\
& &(3^{53},107,53),(5^{18},133,33),(41^{81},163,81),(3^{144},323,161),(5^{249},499,249). 
\end{eqnarray*}
\end{corollary}

\section{Remarks on association schemes}
The results on srgs obtained in Section~\ref{sec4} have implications on association schemes. 
Let $X$ be a finite set. A (symmetric) {\it association scheme} with $d$ classes 
on $X$ consists of sets (binary relations) $R_0,R_1,\ldots,R_d$ which 
partition $X\times X$ and satisfy 
\begin{itemize}
\item[(1)] $R_0=\{(x,x)\,|\,x\in X\}$;
\item[(2)] $R_i$ is symmetric for all $i$;  
\item[(3)] for any $i,j,k\in \{0,1,\ldots,d\}$ there is an integer 
$p_{i,j}^k$ such that given any pair $(x,y)\in R_k$ 
\[
|\{z\in X\,|\,(x,z)\in R_i,(z,y)\in R_j\}|=p_{i,j}^k. 
\] 
\end{itemize}
Note that each of the symmetric relations $R_i$ can be viewed as an undirected 
graph $G_i=(X,R_i)$. Then, the graphs $G_i$, $1\le i\le d$, decompose the 
complete graph with vertex set $X$. An srg with vertex set $X$ and 
its complement form an association scheme on $X$ with two classes. 
If $p_{i,j}^k=p_{j,i}^k$ for all $i,j,k$, then the association scheme is said to be {\it commutative}. 

Let $(X,\{R_i\}_{i=0}^d)$ be a commutative association scheme. For each $i$, $0\leq i\leq d$, let $A_i$ denote the adjacency matrix of $G_i=(X,R_i)$. Then $A_iA_j=\sum_{k=0}^dp_{i,j}^k A_k$ and $A_iA_j=A_jA_i$, for all $0\leq i, j\leq d$.  It follows that $A_0,A_1,\ldots,A_d$ generate a commutative algebra (over the reals) of dimension $d+1$, which is called the 
{\it Bose-Mesner algebra} of the scheme $(X,\{R_i\}_{i=0}^d)$. 
The Bose-Mesner algebra has a unique set of primitive idempotents 
$E_0=(1/|X|)J$, $E_1,\ldots,E_d$, where $J$ is the all-ones matrix. Thus, the algebra has two basis, $\{A_i\,|\,0\le i\le d\}$ and 
$\{E_i\,|\,0\le i\le d\}$. We denote by $P$ the base-change matrix such that $$(A_0,A_1,\ldots,A_d)=(E_0,E_1,\ldots,E_d) \cdot P.$$ The entries in the $i$th column of $P$ are the eigenvalues of $A_i$, $0\leq i\leq d$. The matrix $P$ is called the {\it first eigenmatrix} (or {\it character table}) of the association scheme. 

Given a $d$-class commutative association scheme $(X, \{R_i\}_{0\leq
i\leq d})$, we can take union of classes to form graphs with
larger edge sets (this process is called a {\em{fusion}}). It is
not necessarily guaranteed that the fused collection of graphs will again 
form an association scheme on $X$. If an association scheme has the
property that any of its fusions is also an association scheme, then
we call the association scheme {\em{amorphic}}. A well-known and
important example of amorphic association schemes is given by the
cyclotomic association schemes on $\F_q$ when the cyclotomy is uniform \cite{BMW82}.  For a
partition $\Lambda_0=\{0\}$, $\Lambda_1,\ldots,\Lambda_{d'}\subseteq \{1,2,\ldots,d\}$, let $R_{\Lambda_i}=\bigcup_{k\in \Lambda_i}R_k$. The following 
simple criterion, called {\it the Bannai-Muzychuk criterion}, is very useful for deciding whether $(X,\{R_{\Lambda_i}\}_{i=0}^{d'})$ forms an association scheme or not.
Let $P$ be the first eigenmatrix of the association scheme $(X,\{R_i\}_{i=0}^{d})$. Then, 
$(X,\{R_{\Lambda_i}\}_{i=0}^{d'})$ forms an association scheme if and only if there exists a partition $\Delta_i$, $0\le i\le d'$, of $\{0,1,\ldots ,d\}$, with $\Delta_0=\{0\}$ 
such that each $(\Delta_i,\Lambda_j)$-block of $P$ has a constant row
sum. Moreover, the constant row sum of the $(\Delta_i,\Lambda_j)$-block is the $(i,j)$ entry of the first eigenmatrix of the fusion scheme. 
(For a proof, see \cite{B}.)

From now on, we use the same notation and assumptions as those in 
Lemma~\ref{essen}. 
Let $G_i=\Cay(V,T_i)$, where $T_i=\{i,-i\}\in V^\ast/\{1,-1\}$, and $R_0=\{(x,x)\,|\,x\in V\}$, $R_i=E(\Cay(V,T_i))$. Then, it is obvious that $(V,\{R_i\}_{i=0}^{|V^\ast/\{1,-1\}|})$ forms a 
commutative association scheme. Let $P$ be the first eigenmatrix of this scheme. 
The $(i,j)$ entry of the principal part of $P$ (the matrix obtained by removing the first row and column from $P$) 
of the scheme is given by $\chi_{i}(T_j)$, where both the rows and columns are labeled by the elements of $V^\ast/\{1,-1\}$.   
Write $E_1=D_0\setminus \{0\}$ and $E_{s+2}=D_{C_s^{(e,q)}}$, where $D_0$ and $D_{C_s^{(e,q)}}$
are defined by $D_0=\{x\in V\,|\,Q(x)=0\}$ and $D_{C_s^{(e,q)}}=\{x\in V\,|\,Q(x)\in C_s^{(e,q)}\}$. 
Since $E_i=-E_i$ for $1\le i\le e+1$, 
the subsets  $(\Lambda_{i}=\Delta_i:=)E_i/\{1,-1\}\subseteq V^\ast/\{1,-1\}$, $1\le i\le e+1$, are well defined. 
Then, the row sums of 
the $(\Delta_{i},\Lambda_{j})$-block are given by 
$\chi_{b}(E_j)$, $b\in E_i$. On the other hand, 
Lemma~\ref{essen} implies that for each pair $i,j$, the sum $\chi_{b}(E_j)$ are constant for all 
$b\in E_i$. Thus, by the Bannai-Muzychuk criterion,  the partition gives 
a fusion scheme of $(V,\{R_i\}_{i=0}^{|V^\ast/\{1,-1\}|})$. In summary, we have the following result.

\begin{theorem}\label{asso_semi}
Let $q=p^f$ be a prime power and $n=2m$ be an even positive integer. Let $Q: V=\F_q^n\rightarrow \F_q$ be a nonsingular quadratic form. For any $e\,|\,(q-1)$, let $C_i^{(e,q)}=\omega^i\langle\omega^e\rangle$, $0\leq i\leq e-1$, denote the cyclotomic classes of order $e$ of $\F_q$. Then the decomposition of the 
complete graph on $V$ by $\Cay(V,D_{0}\setminus \{0\})$ and $\Cay(V,D_{C_i^{(e,q)}})$, $0\le i\le e-1$, gives a $(e+1)$-class association scheme.  
\end{theorem}
Next, we give a general sufficient condition for a fusion of the association scheme in Theorem~\ref{asso_semi} to be an association scheme.

\begin{theorem}\label{suff}
Let $q=p^f$ be a prime power and $n=2m$ be an even positive integer. Let $Q: V=\F_q^n\rightarrow \F_q$ be a nonsingular quadratic form. For any $e\,|\,(q-1)$, let $C_i^{(e,q)}=\omega^i\langle\omega^e\rangle$ and $C_{i}^{(e,q^2)}=\gamma^i\langle \gamma^e\rangle$, $0\leq i\leq e-1$, denote the cyclotomic classes of order $e$ of $\F_q$. Assume that  there exists a partition $A_i$, $1\le i\le d$, of $\{i\,|\,0\le i \le e-1\}$ such that  
the decomposition of the complete graph of $\F_{q^2}$ by $\Cay(\F_q,\bigcup_{\ell\in A_i} C_{\ell}^{(e,q^2)})$, $1\le i\le d$, is a fusion scheme of the $e$-class cyclotomic scheme on $\F_{q^2}$. Then, the decomposition of the 
complete graph on $V$ by $\Cay(V,D_{0}\setminus \{0\})$ and $\Cay(V,D_{\bigcup_{\ell\in A_i}C_\ell^{(e,q)}})$, $1\le i\le d$, gives a $(d+1)$-class association scheme.  
\end{theorem}

The above theorem follows immediately from Lemma~\ref{essen}. This can be seen as follows.  By the assumption that the graph decomposition by $\Cay(\F_{q^2},\bigcup_{\ell\in A_i} C_{\ell}^{(e,q^2)})$, $1\le i\le d$, gives a fusion scheme of the $e$-class cyclotomic scheme on $\F_{q^2}$, there exists a partition $\Lambda_h$, $1\le h\le d$, of $\{i\,|\,0\le i \le e-1\}$ such that  for each $1\le h,i\le d$,  $\psi_1'(\gamma^s \bigcup_{\ell\in A_i} C_{\ell}^{(e,q^2)})$ are constant 
for all $s\in \Lambda_h$, i.e.,  
$\phi_b(D_{\bigcup_{\ell\in A_i}C_i^{(e,q)}})$ are constant for all $b\in V$ such that  $Q(b)\in \bigcup_{\ell\in \Lambda_h}C_{\ell}^{(e,q)}$  by Lemma~\ref{essen}. Similarly, $\phi_b(D_{\bigcup_{\ell\in A_i}C_i^{(e,q)}})$ 
 are constant for all $b\in V$ such that  $Q(b)=0$. 
Furthermore, $\phi_b(D_{0})$ is determined according to $Q(b)=0$ or not. 
Thus, by the Bannai-Muzychuk criterion, the conclusion of Theorem~\ref{suff} follows.  We also remark that if the assumed association scheme of $\F_{q^2}$ is amorphic, then so is the resulting scheme on $V$.   

The condition of the above theorem is trivially satisfied in the following case. Let $p$ be a prime, $e>2$, $q=p^{2jr}$, $r\ge 2$, $e\,|\,(p^j+1)$, 
and $j$ is the smallest such positive integer. In this case, since the $e$-class cyclotomic association scheme on $\F_{q^2}$ is amorphic, any fusion of the scheme in Theorem~\ref{asso_semi} forms an association scheme.  
This recovers  Corollary~2.4 of \cite{FWXY}. Also, quite recently, an infinite family of (primitive and non-amorphic) three-class association schemes on $\F_{2^{6s}}$ satisfying the assumption of Theorem~\ref{suff} was found  \cite[Theorem~7 (i)]{FM}.   

Finally, the following theorem of Van Dam~\cite{vD} allows us to put the result of Theorem~\ref{theorem:main2} 
in Section~\ref{sec4} in the context of association schemes. 
\begin{theorem}\label{vD}
Let $\{G_1,G_2,\ldots,G_d\}$ be a decomposition of the complete graph on a set $X$, where each $G_i$ is strongly regular. If $G_i$ are all of Latin square type or all of negative Latin square type, then the decomposition is a $d$-class amorphic association scheme on $X$.  
\end{theorem}
By using Theorem~\ref{theorem:main2} and part (i) of Remark~\ref{re:sec4} in conjunction
with Theorem~\ref{vD}, we have the following: 
\begin{corollary}
Under the same assumptions as in Theorem~\ref{theorem:main2}, the strongly regular decomposition  
$\Cay(V,D_{E})$, $\Cay(V,D_{\F_q^\ast\setminus E})$, $\Cay(V,D_{0}\setminus \{0\})$ yields a $3$-class 
amorphic association scheme. 
\end{corollary}
\section*{Acknowledgments}
The work of K. Momihara was supported by JSPS under Grant-in-Aid for Research Activity Start-up 23840032. The work of Q.  Xiang was done while he is a Program Officer at NSF. The views expressed here are not necessarily those of the NSF.

\end{document}